\documentclass{amsart}

\usepackage{xypic}
\input xy
\xyoption{all}
\usepackage{epsfig}
\usepackage{amsthm}
\usepackage{amssymb}
\usepackage{amsmath}
\usepackage{amscd}

\newcommand{\labell}[1] {\label{#1}}

\newtheorem {Theorem}   {Theorem}

\numberwithin{Theorem}{section}

\newtheorem {Lemma}[Theorem]    {Lemma}
\newtheorem {Proposition}[Theorem]{Proposition}
\theoremstyle{definition}
\newtheorem{Definition}[Theorem]{Definition}
\theoremstyle{remark}
\newtheorem{Remark}[Theorem]{Remark}
\newtheorem{Example}[Theorem]{Example}

\newtheorem {Corollary}[Theorem]{Corollary}
\newtheorem{Claim}[Theorem]{Claim}

\expandafter\chardef\csname pre amssym.def at\endcsname=\the\catcode`\@
\catcode`\@=11
\def\undefine#1{\let#1\undefined}
\def\newsymbol#1#2#3#4#5{\let\next@\relax
 \ifnum#2=\@ne\let\next@\msafam@\else
 \ifnum#2=\tw@\let\next@\msbfam@\fi\fi
 \mathchardef#1="#3\next@#4#5}
\def\mathhexbox@#1#2#3{\relax
 \ifmmode\mathpalette{}{\m@th\mathchar"#1#2#3}%
 \else\leavevmode\hbox{$\m@th\mathchar"#1#2#3$}\fi}
\def\hexnumber@#1{\ifcase#1 0\or 1\or 2\or 3\or 4\or 5\or 6\or 7\or 8\or
 9\or A\or B\or C\or D\or E\or F\fi}

\font\teneufm=eufm10
\font\seveneufm=eufm7
\font\fiveeufm=eufm5
\newfam\eufmfam
\textfont\eufmfam=\teneufm
\scriptfont\eufmfam=\seveneufm
\scriptscriptfont\eufmfam=\fiveeufm

\catcode`\@=\csname pre amssym.def at\endcsname

\def    \eps    {\epsilon}

\newcommand{\supp}{{\mathit supp}\,}
\newcommand{\point}{{\mathit point}}

\newcommand{\id}{{\mathit id}}

\newcommand{\Aa}{{\mathcal A}}

\newcommand{\Ll}{{\mathcal L}}
\newcommand{\Jj}{{\mathcal J}}

\newcommand{\Hh}{{\mathcal H}}

\newcommand{\Mm}{{\mathcal M}}

\newcommand{\Pp}{{\mathcal P}}

\def    \C      {{\mathbb C}}
\def    \R      {{\mathbb R}}
\def    \Z      {{\mathbb Z}}

\def    \ra     {{\rightarrow}}

\def    \12    {{\frac{1}{2}}}

\def    \codim  {\operatorname{codim}}

\def    \pr     {\operatorname{pr}}

\def    \SH     {\operatorname{SH}}
\def    \HF     {\operatorname{HF}}
\def    \CZ     {\operatorname{CZ}}
\def    \CF     {\operatorname{CF}}

\def    \pr     {\preceq}
\def    \wsh       {\operatorname{c_{SH}}}

\def    \ssminus        {\smallsetminus}

\begin{document}






\title[Periodic orbits near symplectic submanifolds]{Symplectic homology and periodic orbits near symplectic submanifolds }

\author[Kai Cieliebak]{Kai Cieliebak}
\author[Viktor L. Ginzburg]{Viktor L. Ginzburg}
\author[Ely Kerman]{Ely Kerman}
\address{Mathematisches Institut, Ludwig--Maximilians--Universitaet,
Theresienstr. 39, 80333 Muenchen, Germany}
\address{Department of Mathematics, UC Santa Cruz, Santa Cruz, CA 95064, USA}
\address{Department of Mathematics, SUNY at Stony Brook, Stony Brook, NY 11794-3651, USA }
\email{kai@mathematik.uni-muenchen.de; ginzburg@math.ucsc.edu; \newline
ely@math.sunysb.edu}

\date{\today}

\thanks{This work was partially supported by the NSF and by the faculty
research funds of the University of California, Santa Cruz. The third author was fully supported by NSERC}


\bigskip

\begin{abstract}
We show that a small neighborhood of a closed symplectic
submanifold in a geometrically bounded aspherical symplectic manifold has
non-vanishing symplectic homology. As a consequence, we establish the existence
of contractible closed characteristics on any thickening of the boundary
of the neighborhood. When applied to twisted geodesic flows on compact
symplectically aspherical manifolds, this implies
the existence of contractible periodic orbits for a dense set of 
low energy values.
\end{abstract}

\maketitle

\section{Introduction}
The symplectic homology of an open subset of a symplectic manifold
is intimately connected with closed characteristics on the shells 
near the boundary of the subset, \cite{fh}. In this paper we show 
that for a small neighborhood of
a compact symplectic submanifold in a geometrically bounded
symplectically aspherical manifold the symplectic homology does not
vanish. As a consequence, we establish the existence of contractible 
closed characteristics on any thickening of the boundary of such a
neighborhood.

Symplectic homology was originally introduced and studied in the
series of papers \cite{cfh,cfhw,fh,fhw}. Here we follow more
recent treatments from \cite{bps} and, more indirectly, from
\cite{vi}. In particular, we consider the (direct limit)
symplectic homology of bounded open sets in a  geometrically
bounded, symplectically aspherical manifold $(W, \omega)$, and we
prove the following:

\begin{Theorem}
\labell{thm:sh-nonzero0}
 Let $U$ be a sufficiently small neighborhood of a compact symplectic 
submanifold $M$ of $W$. Then $\SH^{[a,\,b)}(U)\neq 0$ for $a<b<0$, when
$|a|$ is sufficiently large and $|b|$ is sufficiently small.
\end{Theorem}

This yields the following existence result along the lines of the Weinstein--Moser theorem.
\begin{Theorem}
\labell{thm:orbits}
Let $H\colon W\to \R$ be a smooth function
which attains an isolated minimum on a compact symplectic submanifold
$M$ (say, $H|_M=0$). Then the levels $\{H=\eps\}$ carry contractible 
periodic orbits for a dense set of
small values $\eps>0$.
\end{Theorem}
As an application we obtain a new existence theorem for twisted geodesic flows.
\begin{Corollary}
\labell{cor:twisted} Let $(M, \sigma)$ be a compact symplectically aspherical
manifold and let $H\colon T^*M\to \R$ be the standard kinetic energy Hamiltonian for a Riemannian metric on $M$, i.e., $H(q,p) = \|p\|^2$. Then the 
Hamiltonian flow of $H$ with respect to the twisted form $\omega = \omega_0 + \pi^* \sigma$ has contractible periodic orbits on the levels $\{H=\eps\}$ for a
dense set of small values $\eps>0$.
\end{Corollary}

Theorem \ref{thm:orbits} strengthens or complements some other
recent existence results. 
In particular, it is shown in \cite{gk2},
with no assumptions on the symplectic manifold $(W,
\omega)$ but when $M$ is a Morse--Bott non-degenerate minimum, 
that there exists a sequence of energy values $\eps_k\to 0$
such that all the levels $\{ H=\eps_k\}$ carry contractible
periodic orbits. 
Furthermore, under suitable additional
conditions, a nontrivial topological lower bound for the number of
periodic orbits on every low energy level was established in
\cite{gk1,ke}.

For twisted geodesic flows, Corollary \ref{cor:twisted} augments a
large family of existence results beyond the applications of the
results in \cite{gk1, gk2, ke} described above. For instance, it
was shown in \cite{mac,po} that there are contractible
periodic orbits on a sequence of energy levels $\{ H=\eps_k\}$
with $\eps_k \to 0$ for any compact manifold $M$ and any closed
form $\sigma \neq 0$ with $\sigma|_{\pi_2(M)}=0$. Also, when $M$
is a torus, there are periodic orbits (not necessarily
contractible, e.g., if $\sigma=0$) on almost all energy levels for
any $\sigma$; see \cite{gk1} and references therein, cf.
\cite{ji}. We refer the reader to the survey \cite{gi:survey}
for a discussion of the results on twisted geodesic flows obtained
prior to 1995.

\begin{Remark}
\labell{rmk:assumptions} In Theorem \ref{thm:sh-nonzero0}, the
assumption that $W$ is geometrically bounded can be replaced by
other hypotheses that ensure that the symplectic homology of $U$
is well defined. In particular, the theorem also holds for
symplectic manifolds with  symplectically convex boundary (see
Remark \ref{rmk:convex}). The choice to consider geometrically
bounded manifolds was motivated by our interest in the application
to twisted geodesic flows.
\end{Remark}

The paper is organized as follows. In Section \ref{sec:2}, we
recall the definition of geometrically bounded manifolds and show
that twisted cotangent bundles are geometrically bounded. In
Section \ref{sec:sh}, we briefly review the definitions and
properties of Floer and symplectic homology which are essential
for the proofs of the main theorems. We also show that the
non-vanishing property of symplectic homology implies the
``nearby'' existence of periodic orbits. In Section
\ref{sec:calc} we show that the symplectic homology in question is
non-vanishing and we complete the proofs of the main theorems. Finally,
 in Section \ref{sec:capacity} we discuss the notion of a relative 
homological capacity related to these results.

\subsection*{Acknowledgments.} The authors are deeply grateful to
Ba\c sak G\"urel and Jean-Claude Sikorav for useful discussions
and suggestions. The second author would also like to thank the ESI, Vienna
for the hospitality during the period when this work was essentially
completed. The third author would like to express his
gratitude to the Fields Institute and the University of Toronto
for their hospitality.

\section{Geometrically bounded symplectic manifolds}
\labell{sec:2}

\subsection{Definition of geometrically bounded symplectic manifolds}
\labell{sec:gb} In what follows, we will need to use the Floer
homology of compactly supported functions on non-compact
symplectically aspherical manifolds. In order for this to be
well defined, certain conditions on the manifold at infinity must
be imposed to ensure that the necessary compactness theorems hold.
One standard way to achieve this is to require the symplectic
manifold to be in a certain sense convex. However, the
manifolds we are most interested in -- twisted cotangent bundles
-- do not generally meet this requirement. Hence convexity is
replaced here by the less restrictive requirement that the
manifold is geometrically bounded. Let us recall the definition.

\begin{Definition}
\labell{def:gb}
A symplectic manifold $(W,\omega)$ is said to be
\emph{geometrically bounded} if $W$ admits an almost complex structure $J$
and a complete Riemannian metric $g$ such that
\begin{enumerate}
\item[GB1.] $J$ is uniformly $\omega$-tame, i.e.,
for some positive constants $c_1$ and $c_2$ we have
$$
\omega(X, JX)\geq c_1\| X\|^2 \quad\text{and}\quad
|\omega(X, Y)|\leq c_2\| X\|\,\|Y\|
$$
for all tangent vectors $X$ and $Y$ to $W$.

\item[GB2.] The sectional curvature of $(W,g)$ is bounded from above and
the injectivity radius of $(W,g)$ is bounded away from zero.
\end{enumerate}
\end{Definition}

Observe that if $W$ is compact then condition (GB1) holds
automatically when $g$, $\omega$ and $J$ are compatible, i.e.,
\begin{equation}
\labell{eq:comp}
\omega(\cdot,J \cdot)=g(\cdot, \cdot).
\end{equation}
Clearly, for compact manifolds condition (GB2) also holds for any
metric, and so every compact symplectic manifold is geometrically
bounded. Other examples of geometrically bounded manifolds are
given in the next section. We refer the reader to Chapters V (by
J.-C. Sikorav) and X (by M. Audin, F. Lalonde and L. Polterovich)
in \cite{al} for a more detailed discussion of this concept.

\subsection{Example: twisted cotangent bundles}
\labell{sec:ex-gb} Let us recall the definition of a twisted
cotangent bundle. Consider a closed manifold $M$ and  a closed
two-form (magnetic field) $\sigma$ on $M$. Set $W=T^*M$. The form
$\omega=\omega_0+\pi^*\sigma$ is called a twisted symplectic form,
where $\omega_0$ is the standard symplectic form on $T^*M$ and
$\pi\colon T^*M\to M$ is the natural projection.

One reason these forms are of interest is because the Hamiltonian
flow of the standard kinetic energy Hamiltonian $H\colon T^*M\to
\R$, defined by a Riemannian metric on $M$, describes the motion
of a charge on $M$ in the magnetic field $\sigma$. This is also
called a twisted geodesic flow.

\begin{Proposition}
\labell{prop:gb}
A twisted cotangent bundle $(T^*M,\omega)$ is geometrically bounded for
any closed manifold $M$ and any closed two-form $\sigma$ on $M$.
\end{Proposition}

This proposition is, of course, well known.  The fact that $T^*M$ satisfies condition (GB2) for the natural metric induced by a metric on $M$, is stated without a proof in [AL]; see p. 96 and p. 286.
This assertion is used by Lu in \cite{lu} to show that a twisted cotangent bundle (again, with a natural metric) is geometrically bounded. For the sake of completeness, we outline a proof of Proposition 2.2 below (for a different choice of metric on $T^*M$).\footnote{The authors are
grateful to Jean-Claude Sikorav for suggesting to us the idea of this proof.}


\begin{proof}
Let $\varphi_t$ be the flow on $W=T^*M$ formed by fiberwise dilations
by the factor $e^t$. The standard symplectic form $\omega_0$ is
homogeneous of degree one
with respect to the dilations: $\varphi_t^*\omega_0=e^t\,\omega_0$.
Pick a fiberwise convex hypersurface $\Sigma$ in $W$, enclosing the zero
section $M$. Note that $\Sigma$ has contact type for $\omega_0$, but not
necessarily for the twisted form $\omega$. Denote by $U$ the closure of the unbounded
component of the complement to $\Sigma$ in $W$, i.e.,
$U=\cup_{t\geq 0}\varphi_t(\Sigma)$.

On the vector bundle $TW|_{\Sigma}$, pick any fiberwise metric $g$ and a
complex structure which are compatible with $\omega_0$ in the sense of equation
\eqref{eq:comp}. (We also require the
radial vectors to be $g$-orthogonal to $\Sigma$.) Let us extend
these structures to $U$ so that
\begin{equation}
\labell{eq:hom-m}
\varphi_t^* g=e^t g \quad\text{for}\quad t\geq 0,
\end{equation}
i.e., $g$, just as $\omega_0$,  is homogeneous of degree one with respect to
the dilations, and
$$
J\circ (\varphi_t)_*=(\varphi_t)_*\circ J.
$$
Then the metric $g$, the almost complex structure $J$, and $\omega_0$ are
compatible on $U$ and hence the condition (GB1) holds for these structures.
We will show that this condition also holds for $g$, $J$, and $\omega$
for, perhaps, a dilated hypersurface $\Sigma$.

The metric $g$ is obviously complete. Indeed, identifying $U$ with
$\Sigma\times [1,\infty)$, we see that the metric $g$ has the form
$$
g=\frac{dt^2}{t}+\frac{g|_\Sigma}{t}.
$$
It is then clear that the integral curves $\varphi_t(x)$, for
$t>0$ and $x\in\Sigma$, are minimizing geodesics of $g$. Thus, the
distance from $x$ to $\varphi_t(x)$ is equal to $\ln t$ and goes
to $\infty$ as $t\to\infty$. Therefore, every bounded subset of
$U$ is contained in some shell $\Sigma\times [1,t]$ and is, hence,
relatively compact. This is equivalent to completeness.

It follows readily from \eqref{eq:hom-m} that the sectional curvature of $g$
goes to zero as $x\to \infty$ in $U$. As a consequence, the sectional
curvature of $g$ is bounded from above on $U$. Combined with completeness,
this implies that the injectivity radius is bounded away from zero. Thus,
condition (GB2) holds.

It is also easy to verify condition (GB1). Note that
$\pi^*\sigma$ is homogeneous of degree zero for $\varphi_t$, i.e.,
$\varphi_t^*\pi^*\sigma=\pi^*\sigma$. Then, by a straightforward
calculation, it is clear that for any positive constants $c_1<1$ and $c_2>1$,
this condition holds in the smaller set $\varphi_t(U)$, provided
that $t$ is large enough.

To complete the proof it suffices to extend $g$ and $J$ to the bounded
part $W\ssminus\varphi_t(U)$.
\end{proof}

\begin{Remark}
\labell{rmk:convex} Recall that an open symplectic manifold
$(W,\omega_0)$ is said to be \emph{convex at infinity} if there
exists: a hypersurface $\Sigma\subset W$ which separates $W$ into
one set with compact closure and another, $U$, with non-compact
closure; and a flow of symplectic dilations on $U$, $\varphi_t$
(for $t\geq 0$), which is transversal to $\partial U=\Sigma$. (See
\cite{eg}.) An argument similar to the proof above (with the
exception of the step dealing with $\pi^*\sigma$) shows that a
symplectic manifold $(W,\omega_0)$ which is convex at infinity is
also geometrically bounded. This reasoning together with the fact
that $\pi^*\sigma$ is ``small'' compared to $\omega_0$ is the main
point of the proof above.
\end{Remark}


\section{Symplectic homology}
\labell{sec:sh}

\subsection{Floer theory}
In this section we briefly recall the definitions of Floer homology  
and symplectic homology as well as some of
their properties. All the results here are stated without proof. The 
reader interested in a detailed treatment of Floer homology should consult, 
for example, \cite{hz:book,sa} or the 
original sources \cite{fl1,fl2,fl3}. For the
definition and properties of symplectic homology the reader is referred to
\cite{bps,cfh,cfhw,fh,fhw,vi}.\footnote{The definitions of symplectic 
homology vary considerably from paper to
paper. Here we adopt the approach of \cite{bps}.}

\subsubsection{Floer homology for negative actions.}
Let $(W,\omega)$ be a symplectic manifold, possibly open, which is
geometrically bounded (see Definition \ref{def:gb}). Assume, as well, 
that $(W,\omega)$ is symplectically aspherical, i.e.,
$$
\omega|_{\pi_2(M)}=0 \,\,\,\,\,\,\text{and}\,\,\,\,\,\,
c_1(TM)|_{\pi_2(M)}=0.
$$

Denote by $\mathcal{H}$ the space of smooth, compactly supported
functions on $S^1 \times W$. To each $ H \in \mathcal{H}$ we can
associate the time-dependent Hamiltonian vector field $X_H$ which
is defined by the equation $dH = - i_{X_H} \omega$. The
set of contractible periodic orbits of $X_H$ with period equal to one 
is denoted by $\mathcal{P}(H)$.

Let $\Ll W$ be the space of smooth contractible loops in $W$. We
can also associate to each $H \in \Hh$ the action functional
$\Aa_H \colon \Ll W \ra \R$ given by
$$
\Aa_H(x) = - \int_{D^2} \bar{x}^* \omega +
\int_{S^1}H(t,x)\,dt,
$$
where $ \bar{x} \colon D^2 \ra W$ is any map which restricts to
$S^1 = \partial D^2$ as $x$. This functional is well defined since
$\omega|_{\pi_2 (W)}=0$, and the critical points of $\Aa_H$ are
exactly the elements of $ \mathcal{P}(H)$. The set of critical
values of $\Aa_H$ is called the \emph{action spectrum} and we
denote it by
$$
\mathcal{S}(H)=\{ \Aa_H(x) \mid x \in \mathcal{P}(H) \}.
$$

In general terms, the Floer homology of $H$ is the homology of the
(relative) Morse--Smale--Witten complex of  $\Aa_H$ on $\Ll W$.
However, when $W$ is not compact, every point in the complement of
$\supp H$ is a degenerate $1$-periodic orbit (critical point of
$\Aa_H$) with zero action. To avoid this set, we will only
consider the homology generated by the contractible $1$-periodic
orbits with negative action.

More precisely, for a fixed $a\in (-\infty, 0)$, set
$$
\mathcal{P}^a(H)=\{x \in \mathcal{P}(H) \mid \Aa_H(x) <a\}
$$
and assume that $H$ satisfies the following condition:
\begin{description}
\item [$(*^a)$] Every $1$-periodic orbit $x \in \mathcal{P}^a(H)$ is
nondegenerate.
\end{description}
Since $c_1(TM)|_{\pi_2(M)}=0$, the elements of $\mathcal{P}^a(H)$
are graded by the Conley-Zehnder index $\mu_{\CZ}$; see \cite{sa}.
With this, the Floer complex of $H$ for actions less than $a$ is
the graded $\Z_2$-vector space
$$
\CF^a(H) = \bigoplus_{x \in {\Pp^a (H)}} \Z_2 x.
$$

To define the Floer boundary operator, we first fix an almost
complex structure $J_{gb}$ for which $(W, \omega)$ is
geometrically bounded as in Definition \ref{def:gb}. Let $\Jj$ be
the set of smooth $t$-dependent $\omega$-tame almost complex
structures which are $\omega$-compatible near $\supp(H)$ and are equal
to $J_{gb}$ outside some compact set. Each $J \in \Jj$ defines a
positive-definite bilinear form on $\Ll W$. We can then consider
the moduli space $\Mm (x,y,H,J)$ of \emph{downward} gradient-like
trajectories of $\Aa_H$ which go from $x$ to $y$ and have finite
energy. For a dense subset, $\Jj_{reg}(H) \subset \Jj$, each moduli space
$\Mm (x,y,H,J)$ is a smooth manifold of dimension $\mu_{\CZ}(x)
-\mu_{\CZ}(y)$.

The Floer boundary operator is then defined by
$$
\partial^{H,J} x = \sum_{y \in \mathcal{P}^a(H) \text{ with }
\mu_{\CZ}(x) -\mu_{\CZ}(y) =1} \tau(x,y)y,
$$
where $\tau(x,y)$ stands for the number (mod $2$) of elements in $\Mm
(x,y,H,J)/ \R$ and $\R$ acts (freely) by translation on the
gradient-like trajectories. The operator $\partial^{H,J}$ 
satisfies $\partial^{H,J} \circ
\partial^{H,J} =0$ and the resulting Floer homology groups
$\HF^a(H)$ are independent of the choice of $J \in \Jj_{reg} (H)$.

\begin{Remark}
\labell{rem:compact1} Since $(W, \omega)$ with $J_{gb}$ is
geometrically bounded and $H$ is compactly supported, there is a
uniform $C^0$-bound for the elements of $\Mm (x,y,H,J)$ (see, for
example, Chapter V in \cite{al}). Hence, the compactness of the
appropriate moduli spaces follows from the usual arguments.
\end{Remark}

\begin{Remark}
In this definition of $\partial^{H,J}$ we ignore matters of
orientation by considering only coefficients in $\Z_2$.
\end{Remark}

\begin{Remark}
It is unclear whether $\HF^a(H)$ depends on the choice of
$J_{gb}$. It is independent of this choice if the set of almost
complex structures for which $W$ is geometrically bounded is
connected.
\end{Remark}

It will also be useful to consider Floer homology restricted to
smaller negative action intervals. More precisely, for constants
$-\infty \leq a < b < 0$ let $H \in \mathcal{H}^{a,b}$, where
$$
\mathcal{H} ^{a,b} = \{ H\in \mathcal{H} \mid a,b \notin \mathcal{S} (H) \}.
$$
Assume, as well,  that $H$ has property {\bf $(*^b)$}. Then
$\CF^a(H)$ is a subcomplex of $\CF^b(H)$ and $\HF^{[a,\,b)}(H)$ is
the homology of the quotient complex $\CF^{[a,\,b)}(H) = \CF^b
(H)/ \CF^a (H)$ with the induced boundary operator.

In fact, the set $\mathcal{H}^{a,b}$ is open in $\mathcal{H}$ with
respect to the strong Whitney $C^{\infty}$-topology. Moreover, in
each component of $\mathcal{H}^{a,b}$, the functions with
property {\bf $(*^b)$} form a dense set and have identical Floer homology
groups. Therefore, we can define $\HF^{[a,\,b)}(H)$ for any $H \in
\mathcal{H}^{a,b}$ as the restricted Floer homology of a nearby
function in $\Hh^{a,b}$ with property {\bf $(*^b)$}.

\subsubsection{ Morse--Bott Floer homology}

The extension of the definition of $\HF^{[a,\,b)}(H)$ to every $H \in
\mathcal{H}^{a,b}$ is particularly useful in the Morse--Bott case
which we now describe.

A subset $P \subset \mathcal{P}(H)$ is said to be a
\emph{Morse--Bott manifold of periodic orbits} if the set $C_0 =
\{x(0) \mid x \in P\}$ is a compact submanifold of $W$ and
$T_{x_0}C_0 = \ker (D \phi_H^1(x_0) - \id)$ for every $x_0 \in
C_0.$ Here $\phi_H^1$ is the time-1 flow of $X_H$.

For such sets of periodic orbits we have the following result
which holds for geometrically bounded, symplectically aspherical 
manifolds.

\begin{Theorem}
\labell{thm:mb} {\rm (Po\'zniak, \cite[Corollary 3.5.4]{po};
Biran--Polterovich--Salamon, \cite[Section 5.2]{bps})} Let
$-\infty \leq a < b < 0$ and $H \in \mathcal{H}^{a,b}$. Suppose
that the set $P = \{ x \in \Pp (H) \, \mid \, a < \Aa_H < b\}$ is
a connected Morse--Bott manifold of periodic orbits. Then
$\HF^{[a,\,b)}(H)$ is isomorphic to  $H_* (P; \Z_2)$.
\end{Theorem}

\subsubsection{Monotone homotopies}
\labell{sec:mon-hom}
Let $H,K \in \mathcal{H}^{a,b}$ be two functions with $H(t,x) \geq
K(t,x)$ for all $(t,x) \in S^1 \times W$. Then there exists a monotone
homotopy $s \mapsto K_s$ from $H$ to $K$, i.e., a family of functions $K_s$
such that
$$
K_s=\begin{cases}
H &\text{for $s\in (-\infty,-1]$}\\
K & \text{for $s\in [ 1, \infty)$}
\end{cases}
$$
and $\partial_s K_s \leq 0$. For each such homotopy there is a well-defined
 Floer chain map
$$
\sigma_{KH} \colon \CF^b (H) \ra \CF^b (K).
$$ 
In fact, the map $\sigma_{KH}$ takes $\CF^a(H)$ to  $\CF^a(K)$, and so it 
defines a chain map for the quotient complexes, 
$$
\sigma_{KH} \colon \CF^{[a,\,b)}(H) \ra \CF^{[a,\,b)}(K).
$$ 
This induces a homomorphism of Floer homology,
$$
\sigma_{KH} \colon \HF^{[a,\,b)}(H) \ra \HF^{[a,\,b)}(K).
$$

The following results concerning these homomorphisms are well known ;
see, e.g., \cite{cfh,fh} and \cite[Sections 4.4 and 4.5]{bps}.

\begin{Lemma}
\labell{mon} The homomorphism $\sigma_{KH}$ is independent of the
choice of the monotone homotopy $K_s$ and satisfies the following
identities
\begin{align*}
\sigma_{KH} \circ \sigma_{HG}& = \sigma_{KG} \text{ for } G \geq H \geq K,\\
\sigma_{HH} &= \id  \text{ for every $H \in  \mathcal{H}^{a,b}.$}
\end{align*}
\end{Lemma}

\begin{Lemma}
\label{mo} If $K_s \in \mathcal{H}^{a,b}$ for all $s \in [0,1]$,
then $\sigma_{KH}$ is an isomorphism.
\end{Lemma}

This last result states that the only way in which the map
$\sigma_{KH}$ can fail to be an isomorphism is if periodic orbits,
with action equal to $a$ or $b$, are created during the homotopy.
In fact, this is a particular instance of the following more
general phenomenon.

Let $K_s$ be a monotone homotopy from $H$ to $K$, as above. Assume, for some $c<a$,
that $K \in \mathcal{H}^{c,b}$ and 
$$
\sigma_{KH}(\CF^a(H)) \subset \CF^c(K).
$$ 
Then $\sigma_{KH}$ induces a homomorphism $\hat{\sigma}_{KH}$ from 
$\HF^{[a,\,b)}(H)$ to $\HF^{[c,\,b)}(K)$. (We use the `` $\hat{\,\,}$ '' 
to denote the fact that the map goes between Floer homology groups restricted
to different intervals of actions.) In this case, we get the following 
generalization of Lemma \ref{mo}.

\begin{Lemma}\labell{hat} 
For $s \in [-1,1]$, let $a_s$ be a continuous family of numbers,
less than $b$, such that $a_{-1}=a$ and $a_{1}=c$.
If $K_s \in \mathcal{H}^{a_s,b}$ for every $s \in [-1,1]$, then
$$
\hat{\sigma}_{KH}\colon \HF^{[a,\,b)}(H) \ra \HF^{[c,\,b)}(K)
$$
is an isomorphism.
\end{Lemma}

\begin{Remark}
To define the chain maps above and to prove that the homomorphisms 
$\sigma_{KH}$ and $\hat{\sigma}_{KH}$  are independent of the choice 
of the monotone homotopy (as long as  $K_s \in \mathcal{H}^{a_s,b}$), 
one needs to consider moduli spaces which 
are defined using generic parameterized families in $\Jj$. Each such 
family is constant and equal to $J_{gb}$ outside a compact set in $W$.
Hence, as in Remark \ref{rem:compact1}, the necessary compactness
statements follow from the standard arguments.

\end{Remark}

\subsection{Symplectic homology}

\subsubsection{Direct limits and exhausting sequences}

Following closely the discussion in \cite[Sections 4.6 and
4.7]{bps}, we recall here the algebraic constructions needed to
define symplectic homology. Let $(I, \pr)$ be a partially ordered
set which is thought of as a category with precisely one morphism
from $i$ to $j$ whenever $i \pr j$. A \emph{partially ordered
system of $R$-modules over $I$} is a functor from this category to
the category of $R$-modules. This is written as a pair $(G,
\sigma)$ where $G$ assigns to each $i \in I$ an $R$-module $G_i$
and $\sigma$ assigns to each $i \pr j$ an $R$-module homomorphism
$\sigma_{ji} 
\colon G_i \to G_j$ such that $\sigma_{kj} \circ
\sigma_{ji} = \sigma_{ki}$ for $i \pr j \pr k$ and $\sigma_{ii}$
is the identity on $G_i$.

The partially ordered set $(I, \pr)$ is said to be \emph{upwardly
directed} if for each $i,j \in I$ there is an $l \in I$ such that
$i \pr l$ and $j \pr l$. In this case, the functor $(G, \sigma)$
is called a \emph{directed system of $R$-modules} and its
\emph{direct limit} is defined as
$$\lim_{\longrightarrow}G = \{ (i,x) \mid i \in I, x
\in G_i \}/ \sim ,
$$
where $(i,x) \sim (j,y)$ if and only if there exists an $l
\in I$ such that $i \pr l,\, j \pr l$ and
$\sigma_{li}(x)=\sigma_{lj}(y)$.

Note also that for each $i \in I$ there is a natural homomorphism
$$\iota_i \colon G_i \to \lim_{\longrightarrow} G $$ which takes $x$
to the equivalence class $[i,x]$.

It is possible to actually compute a direct limit given an
\emph{exhausting sequence}. This is a sequence $\{i_{\nu}\}_{\nu \in \Z^+}
\subset I$ with the following two properties
\begin{itemize}
\item For every $\nu \in \Z^+$ we have $i_{\nu} \pr i_{\nu+1}$ and
$\sigma_{i_{\nu+1}i_{\nu}} \colon G_{i_{\nu}} \to G_{i_{\nu+1}}$ is an
isomorphism.
\item For every $i \in I$ there exists a $\nu \in \Z^+$ such that $i
\pr i_{\nu}$.
\end{itemize}
\begin{Lemma}
For such a sequence the map $$\iota_{i_{\nu}} \colon G_{i_{\nu}} \to
\lim_{\longrightarrow} G $$ is an isomorphism for all $\nu \in \Z^+$.
\end{Lemma}

\subsubsection{Construction of $\SH^{[a,\,b)}(U)$}

Let $U \subset W$ be a bounded open subset. We are now in a position 
to define an invariant $\SH^{[a,\,b)}(U)$ of $U$,
called the \emph{symplectic homology of $U$ for the
action interval $[a,\,b)$}. This is done by taking the direct limit of
certain Floer homology groups.

Let
$$
\mathcal{H}^{a,b}(U)= \{ H \in \mathcal{H}^{a,b} \mid
\supp(H) \subset S^1\times U,\, H|_{S^1 \times U} \leq0 \}.
$$
This set has a partial order given by 
$$
H \pr K \iff H \geq K,
$$ 
and it is clearly upwardly directed. 
Note that
$\mathcal{H}^{a,b}(U)\neq \emptyset$ for any $-\infty\leq a< b< 0$. 

To each $H \in \mathcal{H}^{a,b}(U)$ we can associate the
$\Z_2$-vector space given by the Floer homology groups
$\HF^{[a,\,b)}(H)$. Furthermore, for $H \pr K$ we have, by Lemma
\ref{mon}, a homomorphism $\sigma _{KH}$ which satisfies the required
identities to make this a directed system.
Passing to the direct limit, we set
$$
\SH^{[a,\,b)}(U)= \lim_{\longrightarrow} \HF^{[a,\,b)}(H).
$$

As, above, for each $H \in \mathcal{H}^{a,b}(U)$ we have the 
natural homomorphism
$$
\iota_H\colon \HF^{[a,\,b)}(H)\to \SH^{[a,\,b)}(U).
$$

\begin{Remark}
The condition that $\supp (H)\subset S^1\times U$ is sometimes replaced
by the weaker requirement $\supp (H)\subset S^1\times \bar{U}$, see
\cite{cfh,fh}, where the bar denotes the closure. This results in a symplectic
homology with properties similar to the ones considered here and in \cite{bps}.
However, it is likely that the different definitions can yield different
symplectic homology groups when $\partial U$ fails to have contact type.
It is also worth noticing that
the symplectic homology defined with $\supp (H)\subset S^1\times \bar{U}$
cannot be calculated by passing to exhausting sequences (but rather
exhausting families of functions). For
example, the functions $H_\nu$ from Section \ref{sec:near-ex} could not be
used to calculate such symplectic homology.
\end{Remark}

We now recall the definitions and properties of some of the natural
maps introduced in \cite{fh} for symplectic homology. The reader
interested in further details and the proofs of these properties 
should consult \cite{bps,cfh,fh}.

\subsubsection{Monotonicity maps}
Let $V \subset U$ also be a bounded open subset of $W$. Then we have
the natural inclusion
$$
\mathcal{H}^{a,b}(V) \subset \mathcal{H}^{a,b}(U)
$$
and the induced natural homomorphism in
symplectic homology
$$
\phi_{UV}\colon \SH^{[a,\,b)}(V) \to \SH^{[a,\,b)}(U),
$$
which is called a \emph{monotonicity map}.

\begin{Lemma}
\labell{nest}
For a nested sequence $U_1 \subset U_2 \subset U_3$ of bounded open
subsets of $W$,  $\phi_{U_3 U_2} \circ \phi_{U_2 U_1} = \phi_{U_3 U_1}$.
\end{Lemma}

The relationship of the monotonicity map and the map $\iota_H$ is
described, for $H \in \mathcal{H}^{a,b}(V) $, by the following
commutative diagram
\begin{equation}
\labell{i-nu}
\xymatrix@C=0pt{
     \SH^{[a,\,b)}(V) \ar[rr]^{\phi_{UV}}& & \SH^{[a,\,b)}(U)\\
&          \HF^{[a,\,b)}(H) 
\ar[ul]^{\iota_H} \ar[ur]_{\iota_H} &
}
\end{equation}

\subsubsection{Exact triangles}
\labell{sec:ex-tr}
For $-\infty \leq a \leq b \leq c < 0$ and a function $H$ we have
the short exact sequence of complexes
$$
0\;\ra \;\CF^{[a,\,b)}(H)\; \ra
\;\CF^{[a,\,c)}(H)\; \ra \;\CF^{[b,\,c)}(H)\; \ra\; 0.
$$
If $a,b,c \notin \mathcal{S}(H)$,
this sequence yields the exact homology triangle
$\bigtriangleup_{a,b,c}(H)$ given by
$$
\xymatrix@C=0pt{
     \HF^{[a,\,b)}(H) \ar[rr]& & \HF^{[a,\,c)}(H)\ar[dl]^{\Pi} \\
&          \HF^{[b,\,c)}(H) \ar[ul]^{\partial^*}  &
}
$$
Passing to the direct limit, we obtain an exact triangle
$\bigtriangleup_{a,b,c}(U)$ for the symplectic homology for $U$:
$$
\xymatrix@C=1pt{
     \SH^{[a,\,b)}(U) \ar[rr]& & \SH^{[a,\,c)}(U)\ar[dl]^{\Pi} \\
&          \SH^{[b,\,c)}(U) \ar[ul]^{\partial^*}  &
}
$$
Note that the map $\partial^*$ has degree $-1$ in both diagrams. 

Finally, let us describe a relation between the exact triangles and the
homomorphisms induced by monotone homotopies which will be crucial to
us later on. Let $\{K_s\}$ be a monotone homotopy between functions
$H,K \in \mathcal{H}^{a,b}(U)$ with $H \geq K$. As mentioned above, this
induces a homomorphism 
$$ \sigma_{KH} \colon \HF^{[a,b)}(H)
\ra \HF^{[a,b)}(K).$$ 
If, for some $c<a$, we have  $K \in \Hh^{c,b}$ and $\sigma(\CF^a(H)) 
\subset \CF^c(K)$, then
$K_s$ also induces a homomorphism
$$
\hat{\sigma}_{KH} \colon \HF^{[a,\,b)}(H) \ra \HF^{[c,\,b)}(K).
$$
 
\begin{Lemma} The following diagram commutes
\labell{cru}
\begin{equation*}
\xymatrix{
& \HF^{[c,\,b)}(K)  \ar[d]^{\Pi}\\
\HF^{[a,\,b)}(H) \ar[ur]^{\hat{\sigma}_{KH}} \ar[r]^{\sigma_{KH}}
& \HF^{[a,\,b)}(K)
 }
\end{equation*}

where $\Pi$ is the map from the exact triangle
$\bigtriangleup_{c,a,b}(K)$.
\end{Lemma}

\begin{proof}
Using the explicit definition of the maps $\sigma_{KH}$ and
$\hat{\sigma}_{KH}$ (see \cite{cfh,fh}),
it is not hard to see that the diagram commutes even on the chain level.
\end{proof}

\subsection{Symplectic homology and ``nearby existence'' of periodic orbits}
\labell{sec:near-ex}
 
Let $U$ be a bounded open subset of $W$ with smooth boundary
$\partial U$. Then the symplectic homology $\SH^{[a,\,b)}(U)$ contains 
information about the symplectic geometry of the hypersurface $\partial U$. 
In particular, it detects the contractible closed characteristics on or
near $\partial U$.

To make this idea more precise, set $\Sigma = \partial U$ and define a
\emph{thickening of $\Sigma$} to be a map $\Psi \colon
[-1,0] \times \Sigma \ra M $ which is a diffeomorphism on its image and
satisfies $\Psi( \{0\} \times \Sigma ) =
\Sigma$ and $\Psi( [-1,0) \times \Sigma ) \subset U$. The following result, 
in a slightly different version, essentially
goes back to \cite{fh}.

\begin{Proposition}
\labell{near}
Suppose that $\SH^{[a,\,b)}(U) \not= 0$ for
$-\infty < a \leq b \leq 0$.
Then for any thickening $\Psi$ of $\Sigma = \partial U$ there is a
sequence $\{t_i\} \subset [-1,0)$, converging to $0$, such that the
corresponding hypersurfaces $\Sigma_{t_i}=\Psi({t_i}\times \Sigma)$ all contain contractible closed characteristics.
\end{Proposition}

\begin{proof}
Assume that no such sequence exists. Then there exists a number
$\lambda \in [-1,0)$ such that none of the hypersurfaces $\Sigma_t$
for $t \in [\lambda,0)$ contains a contractible closed
characteristic. Consider a sequence of smooth functions $f_{\nu}
\colon [\lambda, 0] \ra (- \infty, 0]$ with the following properties
\begin{itemize}
\item $f_{\nu}(t)=\begin{cases}
0 &\text{for $t\in [\lambda/4{\nu},0]$},\\
a-{\nu} & \text{for $t\in [ \lambda, \lambda(1-1/4{\nu})]$.}
\end{cases}
$ \\
\item $f'_{\nu}(t)>0 \quad \text{for}\quad t \in
(\lambda(1-1/4{\nu}),\lambda/4{\nu}).$\\
\item $f_{\nu}(t) \geq f_{\nu+1}(t)$ for all $\nu \in \Z^+$ and $t
\in [\lambda, 0]$.
\end{itemize}

Using these functions, we construct a sequence of Hamiltonians
$H_{\nu} \colon M \ra \R$ by setting
$$
H_{\nu}(x)=\begin{cases}
0 & \text{for $x \in  M \setminus \bar{U}$},\\
f_{\nu}(t) & \text{for $x \in \Sigma_t, t \in [-1,0)$}, \\
a-{\nu} & \text{for $x\in U \setminus \Psi([-1,0)\times \Sigma)$}.
\end{cases}
$$
The functions $\{H_{\nu}\}$ clearly have the second property of an
exhausting sequence and satisfy $H_{\nu} \geq H_{\nu+1}$. We claim
that the maps $\sigma_{H_{\nu+1}H_{\nu}}$ are (trivially)
isomorphisms and the functions $\{H_{\nu}\}$ indeed form an
exhausting sequence. By construction, the Hamiltonian vector
fields $X_{H_{\nu}}$ are only 
nonzero 
on the level sets corresponding
to the hypersurfaces $\{\Sigma_t\}_{t \in [\lambda, 0)}$. On these
level sets, the trajectories of the $X_{H_{\nu}}$ correspond to
characteristics and, by assumption, none of these can be closed
and contractible. Away from the level sets $\{\Sigma_t\}_{t \in
[\lambda, 0)}$, the trajectories are all constant and have action
either equal to zero or less that $a$. Thus, the Floer homology
groups $\HF^{[a,\,b)}(H_{\nu})$ are trivial for all ${\nu}$ and
the functions $H_{\nu}$ form an exhausting sequence. This implies that
$$
\SH^{[a,\,b)}(U) = \iota_{H_{\nu}}(\HF^{[a,\,b)}(H_{\nu}))=0,
$$
which is a contradiction.
\end{proof}

\section{Proofs of the main theorems}
\labell{sec:calc} 
We are now in a position to prove Theorem
\ref{thm:sh-nonzero0} and Theorem \ref{thm:orbits}. Since most of
this section is devoted to the proof of the first of these
theorems, we recall its assertion for the sake of convenience.

\begin{Theorem}
\labell{thm:sh-nonzero} Let $M$ be a compact symplectic
submanifold of $(W, \omega)$, a geometrically bounded symplectically 
aspherical manifold. Let $U$ be a sufficiently small neighborhood
of $M$ in $W$. Then $\SH^{[a,\,b)}(U)\neq 0$ for all sufficiently
large negative $a$ and small negative $b$.
\end{Theorem}

Combining Theorem \ref{thm:sh-nonzero} and Proposition \ref{near},
we immediately obtain Theorem \ref{thm:orbits}, i.e., the existence of
periodic orbits near any compact oriented hypersurface bounding a small
neighborhood of $M$ in $W$.

\subsection{Proof of Theorem \ref{thm:sh-nonzero}}
\labell{sec:pf}
Let $E\to M$ be the symplectic normal bundle to $M$, i.e., the symplectic orthogonal
complement to $TM$ in $TW|_M$. We first recall that a neighborhood
of the zero section in $E$ has a natural symplectic structure. Moreover,
on this neighborhood, there exists a fiberwise quadratic Hamiltonian 
whose flow is periodic. This can be seen as follows.

Let us equip $E$ with a Hermitian metric which is compatible with
the fiberwise symplectic structure on $E$. We denote by
$\rho\colon E\to \R$ the square of the fiberwise norm, i.e.,
$\rho(z)=\|z\|^2$. Recall that $E$ has a canonical fiberwise
one-form whose differential is the fiberwise symplectic form. 
(The value of this form at $z\in E$ is equal to the contraction of
the fiberwise symplectic form by $z$.) 

Fixing a Hermitian connection
on $E$ we extend this fiberwise one-form to a genuine one-form
$\theta$ on $E$. Then the form
$$
\omega_E=\12 d\theta+\sigma
$$
is symplectic on a neighborhood of the zero section in $E$. Here we have identified
$\sigma=\omega|_M$ with its pull-back to $E$.

It is not hard to see that all the orbits of the Hamiltonian flow
of the function $ \rho \colon E\to\R$ are periodic with period
$\pi$, just as for the square of the standard norm on $\R^{2n}$.
(In fact, $i_X\omega_E= z/2$ everywhere on $E$, where $X$ is the
vector field generating the standard fiberwise Hopf action.) This
fact will be essential for the calculation of the Floer homology
of a small tubular neighborhood of $M$ in $W$.

By the symplectic neighborhood theorem, a sufficiently small
neighborhood of $M$ in $(W,\omega)$ is symplectomorphic to a small
neighborhood of $M$ in $(E,\omega_E)$. From now on, we will assume
that this identification has been made. Hence, in what follows,
$\omega=\omega_E$ and $\rho$ is regarded as a function on a
neighborhood of $M$ in $W$. Sometimes, we will write $\rho(z)$ as
$\| z \|^2$.

Denote by $B_r$ the disc bundle of radius $r$ in $E$. 
The key to the proof of Theorem \ref{thm:sh-nonzero} is the following result.

\begin{Proposition}
\labell{cont}
Let $0< r < R$ be sufficiently small and assume that
$$
a \in [-\infty, - \pi R^2)\quad\text{and}\quad b \in [- \pi r^2, 0).
$$
Then
$$
\SH^{[a,\,b)}_{n_0}(B_r)= \SH^{[a,\,b)}_{n_0}(B_R) = \Z_2,
$$
where
$$
n_0=\frac{1}{2}(\dim M-\codim M),
$$
and the monotonicity map
$\phi_{B_R B_r} \colon \SH^{[a,\,b)}_{n_0}(B_r) \to \SH^{[a,\,b)}_{n_0}(B_R)$
is an isomorphism.
\end{Proposition}

\begin{Remark}
In \cite{fhw}, Floer, Hofer and Wysocki compute the symplectic
homology of ellipsoids in $\R^{2n}$ and prove results which are
similar to Proposition \ref{cont} for $E=\R^{2n}$ (see Corollary 2
and Proposition 6). Our proof of Proposition \ref{cont} is
similar in spirit to that in \cite{fhw}. However, instead of
using explicitly perturbed non-degenerate Hamiltonians as in
\cite{fhw}, we use Morse--Bott Floer theory (see Theorem
\ref{thm:mb}) to calculate the relevant Floer homology groups.
This is motivated by the successful use of this technique in
\cite{bps}.
\end{Remark}

\begin{Remark}
The hypothesis that $R>0$ is small is only needed to guarantee that
$B_R$ is contained in $W$ for which the Floer homology is
defined. This hypothesis can sometimes be relaxed. For example,
if $\omega_E$ is symplectic on the entire total space of $E$ and
$(E,\omega)$ is geometrically bounded, we can take arbitrarily large
radii $R$ and $r$.
\end{Remark}

Theorem \ref{thm:sh-nonzero} follows immediately from Proposition
\ref{cont}. Indeed, when $U$ is small enough, we may assume that
$U\subset B_R$ and, since $M\subset U$, the disc bundle $B_r$ is
contained in $U$ for sufficiently small $r>0$, i.e., 
$$
B_r\subset U\subset B_R.
$$
By Lemma \ref{nest}, the inclusion map $\phi_{B_R B_r}$ factors as
$$
\SH^{[a,\,b)}(B_r) \to \SH^{[a,\,b)}(U)\to \SH^{[a,\,b)}(B_R),
$$
and by Proposition \ref{cont}, the map $\phi_{B_R B_r}$ is a nontrivial
isomorphism for suitable negative
constants $a$ and $b$. Thus, we conclude that $\SH^{[a,\,b)}(U)\neq 0$.


\subsection{Proof of  Proposition \ref{cont}}
\subsubsection{Outline of the proof}
First we consider the following diagram which allows us to work in the setting of Floer homology:
\begin{equation}
\labell{sh-hf}
\xymatrix{
\SH^{[a,\,b)}(B_r) \ar[r]^{\phi_{B_R B_r}} &  \SH^{[a,\,b)}(B_R)\\
\HF^{[a,\,b)}(H_{\nu}^r)
\ar[u]^{\iota_{H_{\nu}^r}}
\ar[r]^{\sigma_{H^R_\nu H^r_\nu}} &
\HF^{[a,\,b)}(H_{\nu}^R)
\ar[u]_{\iota_{H_{\nu}^R}}.
}
\end{equation}

Here, $\{H_{\nu}^r\}_{\nu \in \Z^+}$ and $\{H_{\nu}^R\}_{\nu
\in \Z^+}$ are chosen to be upward exhausting sequences for 
the symplectic homology groups $\SH^{[a,\,b)}(B_r)$ and 
$\SH^{[a,\,b)}(B_R)$, respectively. Hence, the two vertical 
arrows are isomorphisms (for sufficiently large $\nu$). There is also 
a natural monotone homotopy $H_{\nu}^s$
from $H_{\nu}^r$ to $H_{\nu}^R$, for each $\nu$,  so the map 
$\sigma_{H_{\nu}^R H_{\nu}^s}$
is well-defined. The fact that the diagram commutes then follows easily from
diagram \eqref{i-nu} and the definition of the direct limit. 

 By diagram \eqref{sh-hf}, Proposition \ref{cont} will follow from the two 
results below which concern the restricted Floer homology
of the functions $H^s_{\nu}$. In particular, these results imply that the 
groups $\HF^{[a,\,b)}_{n_0}(H_{\nu}^r)$ and $\HF^{[a,\,b)}_{n_0}(H_{\nu}^R)$ 
are nontrivial, and that $\sigma_{H^R_\nu H^r_\nu}$ 
is an isomorphism between them.

\begin{Claim}
\label{claim1} Let $n_0=\frac{1}{2}(\dim M-\codim M)$. For every $s\in [r,R]$
 and any $c<a$ which is not a negative integer multiple of $\pi s^2$, the map
$$
\Pi \colon \HF_{n_0}^{[c,\,b)}(H^s_{\nu}) \to
\HF_{n_0}^{[a,\,b)}(H^s_{\nu})
$$
is (eventually) an isomorphism and
$$\HF_{n_0}^{[c,\,b)}(H^s_{\nu}) =
\HF_{n_0}^{[a,\,b)}(H^s_{\nu}) = \Z_2.
$$
\end{Claim}

Here, $\Pi$ is the map from the exact triangle $\bigtriangleup_{c,a,b}(H^s_{\nu})$. 
\begin{Claim}
\label{claim2} Let $k'$ be the largest integer in $(0, -a/ \pi r^2)$ and let
$s_0,s_1 \in [r,R]$ satisfy 
$$
s_0 < s_1 < \sqrt{1+ 1/(k'+1)} s_0.
$$ 
Then there are constants $c$, satisfying $c <a$, such that the map
$$
\hat{\sigma}_{H^{s_1}_{\nu}H^{s_0}_{\nu}} \colon \HF^{[a,\,b)}(H_{\nu}^{s_0})
\to \HF^{[c,\,b)}(H_{\nu}^{s_1})
$$ 
is well-defined and an isomorphism. Moreover, $c$ can be chosen so 
that $c \notin - \pi s_1^2 \Z^+$.
\end{Claim}

Claim \ref{claim1} implies that the groups $\HF^{[a,\,b)}_{n_0}(H_{\nu}^r)$ and $\HF^{[a,\,b)}_{n_0}(H_{\nu}^R)$ are both isomorphic to $\Z_2$. It just remains to show that 
$$
\sigma_{H^{R}_{\nu} H^{r}_{\nu}} \colon \HF^{[a,\,b)}_{n_0}(H_{\nu}^r)
\to \HF^{[a,\,b)}_{n_0}(H_{\nu}^R)
$$ is an isomorphism. To see this, let $s_0$ and $s_1$ satisfy 
the assumptions of Claim \ref{claim2} and consider the following 
version of the commutative diagram from
Lemma \ref{cru}:
\begin{equation}
\labell{tri2}
\xymatrix{
& \HF^{[c,\,b)}(H_{\nu}^{s_1})  \ar[d]^{\Pi}\\
\HF^{[a,\,b)}(H_{\nu}^{s_0})
\ar[ur]^{\hat{\sigma}_{H^{s_1}_{\nu} H^{s_0}_{\nu}}}
\ar[r]^{\sigma_{H^{s_1}_{\nu} H^{s_0}_{\nu}}}
& \HF^{[a,\,b)}(H_{\nu}^{s^1})
}
\end{equation}
Together, Claims \ref{claim1} and \ref{claim2} imply that the map 
$$
\sigma_{H^{s_1}_{\nu} H^{s_0}_{\nu}} \colon \HF_{n_0}^{[a,\,b)}(H_{\nu}^{s_0})
 \ra \HF_{n_0}^{[a,\,b)}(H_{\nu}^{s_1})
$$
is an isomorphism. Now, by Lemma \ref{mon}, for any sequence 
$\{s_i\}_{i=1,\dots ,k}$ with 
$$
r=s_0< s_1< \dots < s_{k-1}<s_k=R,
$$ 
we can decompose $\sigma_{H^{R}_{\nu} H^{r}_{\nu}}$ as
$$
\sigma_{H^{R}_{\nu} H^{r}_{\nu}}= \sigma_{H^{R}_{\nu} H^{s_{k-1}}_{\nu}} \circ
\sigma_{H^{s_{k-1}}_{\nu} H^{s_{k-2}}_{\nu}} \circ \dots \circ \sigma_{H^{s_1}_{\nu} H^{s_{r}}_{\nu}}.
$$
Choosing the $s_i$ so that
$s_i < s_{i+1} < \sqrt{1+1/(k'+1)}s_i$, it follows that
$\sigma_{H^{R}_{\nu} H^{r}_{\nu}}$ is an isomorphism from 
$\HF^{[a,\,b)}_{n_0}(H_{\nu}^r)$ to $\HF^{[a,\,b)}_{n_0}(H_{\nu}^R)$. 

\subsubsection{Rescaling} For simplicity, we will assume that $r=1$ and $R >1$. This can always be achieved by scaling the symplectic structure on $W$ and
the metric on $E$. In this case, the fixed constants $a$ and $b$ satisfy
$a <-\pi R^2 < -\pi < b < 0$, and $k'$ is the largest integer in $(0, -a/\pi)$.

\subsubsection{The functions $H^s_{\nu}$}
\labell{sec:exh}
Here, we construct the functions $H^s_{\nu}$, describe their
closed orbits with period equal to one, and prove that for each 
$s \in [1,\,R]$ the sequence $\{H^s_{\nu}\}$ is (eventually) upwardly 
exhausting for $\SH^{[a,\,b)}(B_s)$.

We begin by constructing the exhausting sequence $\{H_{\nu}\} \equiv
\{H^1_{\nu}\}$ for $\SH^{[a,\,b)}(B_1)$. These functions depend only on the 
norm of the fibre variable, i.e.,  $H_{\nu}=H_{\nu}(\|z\|)$. The 
functions $H^s_\nu$ are then obtained by rescaling the argument as 
$H^s_\nu(\|z\|)= H_\nu(\|z\|/s)$.

Consider a sequence of smooth nondecreasing functions $f_{\nu} \colon
[0, \infty) \to (-\infty, 0]$, for $\nu \in \Z^+$, with the following
properties (see Figure \ref{fig1}):
\begin{enumerate}
\item[(F1)]
$f_{\nu}(t)=\begin{cases}
0 &\text{for $t\in [1-1/2^{\nu+2},\,\infty),$}\\

f_{\nu}(0) & \text{for $t\in [0, \,1-1/2^{\nu}].$}
\end{cases}$\\
\item[(F2)]
$f_{\nu}(t)\in [-1/2^{\nu+2},0]$ for $t \in (1-1/2^{\nu+1},\,1-1/2^{\nu+2}).$\\
\item[(F3)]
 $f'_ {\nu}(t)= 2^{\nu+3}$ for all $t\in [1-3/2^{\nu+2},\,1-1/2^{\nu+1}]$.\\
\item[(F4)]
$f''_ {\nu}> 0$ for $t \in (1-1/2^{\nu},\,1-3/2^{\nu+2})$  and\\
 $f''_{\nu} <  0$ for $t \in (1-1/2^{\nu+1},\,1-1/2^{\nu+2}).$
\end{enumerate}

\begin{figure}[htbp]
\caption{The functions $f_{\nu}$ }
\labell{fig1}
\begin{center}
 \psfig{file=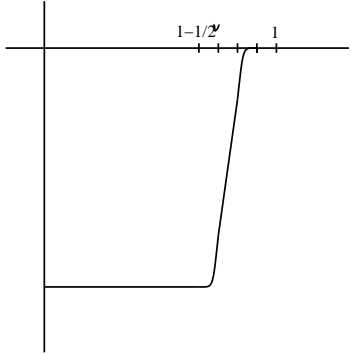}
\end{center}
\end{figure}

On $B=B_1$, set
$$
H_{\nu} (z) = f_{\nu}(\|z\|^2),
$$
and extend this function to be identically zero on $W\ssminus B$.

Recall that on the neighborhood of $M$ where $\omega_E$ is
symplectic, the Hamiltonian flow of $\rho(z)=\|z\|^2$ is totally
periodic with period $\pi$. Hence, for each ${\nu}$, the flow of
$X_{H_{\nu}}$ is totally periodic and the trajectories on each
level set of $H_{\nu}$ all have the same period. The level sets which 
consist of non-constant closed orbits with period equal to one 
correspond to the solutions of the equations
\begin{equation}
\labell{k}
f'_{\nu}(\|z\|^2)= k\pi,
\end{equation}
where $k$ is a positive integer. These solutions are easy to classify. 
When $k$ is in the 
interval $[1, 2^{\nu+3}/ \pi)$, there 
are exactly two levels of $H_\nu$ where equation \eqref{k} holds:
\begin{eqnarray*}
\|z\|^2 &=& t^-_{\nu,k},\quad \text{where
$t^-_{\nu,k} \in [1-1/2^{\nu},1-3/2^{\nu+2})$, and}\\
\|z\|^2 &=& t^+_{\nu,k},\quad \text{where
$t^+_{\nu,k} \in (1-1/2^{\nu+1},1-1/2^{\nu+2}]$.}
\end{eqnarray*}
Since $2^{\nu+3} \notin \pi \Z^+$, there are no other solutions
to \eqref{k}.

\begin{Remark}
\labell{remark:large-k}
As $\nu$ increases, one just gets new solutions, $t^{\pm}_{\nu,k}$, of 
equation \eqref{k} for increasingly large values of $k$.
\end{Remark}

We now show that for large enough $\nu$ that only the orbits on
the levels $\|z\|^2 = t^+_{\nu,k}\in (1-1/2^{\nu+1},\,1-1/2^{\nu+2})$
with $k\in [1,\,-a/\pi)$ have action in $(a,\,b]$. By a simple calculation, 
the action of each of the orbits on the
level set $\|z\|^2=t^{\pm}_{\nu,k}$ is equal to
$$
\Aa^{\pm}_{\nu,k} (1) = f_{\nu}(t^{\pm}_{\nu,k}) - t^{\pm}_ {\nu,k}\pi k.
$$
Since $f_{\nu}(t)< -2^{\nu}$ for all $t \in [0, 1-3/2^{\nu+2}]$, the
actions $\Aa^-_{\nu,k}$ decrease exponentially to negative infinity as 
$\nu \to \infty$. Hence, the orbits on the levels $\|z\|^2=t^-_{\nu,k}$ 
will eventually have action less than $a$. It is also easy to check 
that the actions $\Aa^{+}_{\nu,k} (1)$ decrease monotonically to 
$-\pi k$ as $\nu \to \infty$. Thus, because $b > -\pi$, all the levels
$\|z\|^2=t^+_{\nu,k}$ will eventually have action less than $b$.
Similarly, the levels $\|z\|^2=t^+_{\nu,k}$, with $k \in (-a/\pi,
2^{\nu+3}/\pi)$, will eventually have action less than $a$.

From this discussion we also see that the functions $\{H_{\nu}\}$
belong to $\mathcal{H}^{a,b}$ for all sufficiently large $\nu$.

\begin{Lemma}
\labell{lemma:exh}
The sequence $\{H_{\nu}\}$ is (eventually) upwardly exhausting for the
symplectic homology $\SH^{[a,\,b)}(B)$.
\end{Lemma}

\begin{proof}
The only property that is not immediately obvious is that
$\sigma_{H_{{\nu}+1}H_{\nu}}$ is an isomorphism for all
sufficiently large ${\nu}$. To see this, consider $\nu$ as a
continuous parameter and look at the monotone homotopy
$\{H_{{\nu}+\tau}\}_{\tau \in [0,1]}$ between $H_{\nu}$ and
$H_{\nu+1}$. When $\nu$ is large enough,
the new periodic orbits with period one which appear as $\tau$
goes from zero to one will all have action less than $a$.
(This follows from Remark \ref{remark:large-k} and the fact that 
$\Aa^{+}_{\nu,k} (1) \to -k\pi$ as $\nu \to \infty$.)
Hence, $\{H_{{\nu}+\tau }\} \in \mathcal{H}^{a,b}$ for all $\tau \in
[0,1]$ and Lemma \ref{lemma:exh} immediately follows from Lemma
\ref{mo}.
\end{proof}

\begin{Remark}
The same argument shows that $\{H_{\nu}\}$ is (eventually)
upwardly exhausting for every $\SH^{[c,\,d)}(B)$ with $ a \leq c
\leq d <0$.
\end{Remark}


For $1 \leq s \leq R$, we now define the functions $H^s_{\nu}$ by
the equation
$$
H^s_{\nu}(\|z\|) = H_{\nu}(\|z\|/s) = f_{\nu}(\|z\|^2/s^2).
$$
For these new Hamiltonians the level sets consisting of $1$-periodic
orbits are defined by the equations
\begin{equation}
\label{k,s}
f'_{\nu}(\|z\|^2/s^2)= s^2 k\pi.
\end{equation}
When $k$ is any integer in $[1,2^{\nu+3}/ \pi s^2)$ each of these
equations has two solutions, $\|z\|^2 = s^2 t_{\nu,k}^{\pm} (s)$,
where $ t^{-}_{\nu,k}(s) $ is in $[1-1/2^{\nu},1-3/2^{\nu+2})$ and
$t^+_{\nu,k}(s)$ is in $(1-1/2^{\nu+1},1-1/2^{\nu+2}]. $ We also
have $t_{\nu,k}^{\pm} (s) \to t_{\nu,k}^{\pm}$ as $s \to 1$ for
all $k$ and $\nu$.

The actions of the orbits on these level sets are equal to
$$
\Aa^{\pm}_{\nu,k} (s) = f_{\nu}(t^{\pm}_{\nu,k}(s)) -
t^{\pm}_{\nu,k}(s) \pi k s^2.
$$
As $\nu \to \infty$, the actions  $\Aa^-_{\nu,k}(s)$ decrease exponentially
to negative infinity, and the actions $\Aa^+_{\nu,k}(s)$ decrease monotonically to $-\pi k s^2$.

Just as above, $H^s_{\nu}$ is  in $\Hh^{a,\,b}$ for sufficiently large
$\nu$, and the only orbits with action in $[a,b)$ are those on the
level sets $\|z\|^2=s^2t^+_{\nu,k}(s)$ for $k \in [1, \,-a/ \pi
s^2)$. The same arguments also show that the functions $H^s_{\nu}$
(eventually) form an upward exhausting sequence for
$\SH^{[a,\,b)}(B_s)$. In particular, this is true for $s=R$.

\subsubsection{ Proof of Claim \ref{claim1}} 

By the scaling properties of the functions $H^s_{\nu}$ it suffices to prove 
Claim \ref{claim1} for the case $s=1$. 
Recall that the periodic orbits of $H_{\nu}$ which lie on the level 
sets defined by $\|z\|^2 = t^-_{\nu,k}$ have actions which decrease
exponentially as $\nu \to \infty$. So, when considering bounded intervals 
of actions and large enough $\nu$, the only relevant periodic orbits 
are those on the level sets of $H_{\nu}$ defined by $\|z\|^2 = t^+_{\nu,k}$.
The actions of these orbits is given by
$$
\Aa^+_{\nu,k} (1)=f_{\nu}(t^+_{\nu,k})-t^+_{\nu,k}\pi k
$$
which decreases to $-\pi k$ monotonically. The following result is now obvious.

\begin{Lemma}
\labell{empty} Let $l$ be a non-negative integer and let $-(l+1) \pi <
\alpha \leq \beta \leq -l \pi$. Then for sufficiently large $\nu $
$$
\HF^{[\alpha,\,\beta)} (H_{\nu}) =0.
$$
\end{Lemma}

Now, consider intervals of the form $[-k\pi -\eps, -k\pi + \delta)$ where
$0 < \eps, \delta <\pi$ and $k \in \Z^+$. By the same reasoning as
above, for sufficiently large ${\nu}$, the only level set with action
in this interval is $\|z\|^2=t^+_{\nu,k}$.

\begin{Lemma}
\labell{shell}
For any $0 < \eps, \delta <\pi$ and all $\nu$ sufficiently large
$$
\HF^{[-k\pi -\eps,-k\pi + \delta)}_*(H_{\nu})
 \cong H_{*+j(k)}(SM, \Z_2),
$$
where $SM$ is the sphere bundle of the normal bundle to $M$ in $W$
and
$$
j(k)=\frac{1}{2}\dim W+k\codim M-1.
$$
\end{Lemma}

\begin{proof}
The assertion, up to the value of the shift $j$, follows immediately from
Theorem \ref{thm:mb} provided that the
set of periodic orbits in $\|z\|^2=t^+_{\nu,k}$ is Morse--Bott. However,
it is straightforward to check that the Morse--Bott condition is
equivalent, for the Hamiltonians in questions, to the
condition $f''_{\nu}(t^+_k) \not= 0$, \cite[Lemma 5.3.2]{bps}. The value
of $j(k)$ is determined by the standard analysis of the Conley--Zehnder indices
for a small generic perturbation of the Hamiltonian.
\end{proof}

Proceeding with the proof of Claim \ref{claim1}, we first show
that for any constant $c<a$, where $c$ is not a negative
integer multiple of $\pi$, the map
$$
\Pi\colon \HF_{n_0}^{[c,\,b)}(H_{\nu})\to \HF_{n_0}^{[a,\,b)}(H_{\nu})
$$
is an isomorphism provided that $\nu$ is large enough.

To do this, we consider the following exact sequence which is a
part of the exact triangle $\bigtriangleup_{c,a,b}(H_\nu)$ (see
Section \ref{sec:ex-tr}).
\begin{equation}
\labell{ex}
\HF_{n_0}^{[c,\,a)}(H_{\nu})\to
\HF_{n_0}^{[c,\,b)}(H_{\nu})\stackrel{\Pi}{\to}
\HF_{n_0}^{[a,\,b)}(H_{\nu})\;\to \HF_{n_0-1}^{[c,\,a)}(H_{\nu}).
\end{equation}

First, assume that $c$ and $a$ are contained in
an  interval of the form $(-\pi(k+1),-\pi k]$. Then Lemma
\ref{empty} immediately implies that the first and last terms in
\eqref{ex} vanish for sufficiently large $\nu$ and hence $\Pi$ is an
isomorphism.

Next, we look at the case where $c$ and $a$ are in adjacent
intervals, i.e.,
$$
-\pi(k+1)<c< -\pi k < a \leq -\pi (k-1)
$$
for some integer $k\geq 2$ (since $a< -\pi$). By Lemma \ref{shell}, if $\nu$ is
large enough,
$$
\HF_{n_0}^{[c,\,a)}(H_{\nu})=H_{n_0+j(k)}(SM,\Z_2)=0,
$$
because $n_0+j(1)=\dim SM$ and 
$j(1)<j(k)$. (The sequence $j(k)$ is increasing.) In a
similar manner we obtain
$$
\HF_{n_0-1}^{[c,\,a)}(H_{\nu})=H_{n_0+j(k)-1}(SM,\Z_2)=0,
$$
when $k\geq 2$. Hence, in this case, $\Pi$ is again an isomorphism.

Now, for any $c<a$ such that $c \notin \pi \Z^+$, the map $\Pi$
can be factored as a composition of isomorphisms for points in
adjacent intervals. More precisely, $\Pi = \Pi_{l} \circ \Pi_{l-1}
\circ \cdots \circ \Pi_{1}$ where
$$
\Pi_{j} \colon \HF^{[c_j, b)}_{n_0}(H_{\nu}) \to \HF^{[c_{j-1},
b)}_{n_0}(H_{\nu}).
$$
The constants $c_j$ are chosen to be in adjacent intervals so that
$c_0= a$ and $c_j \notin - \pi \Z^+$ for $j \geq 1$. The previous
argument then implies that each $\Pi_j$ is an isomorphism for
sufficiently large $\nu$. Thus, the same is true for $\Pi$.

It remains to show that
$\HF_{n_0}^{[c,\,b)}(H_{\nu})=\HF_{n_0}^{[a,\,b)}(H_{\nu})=\Z_2$.
The above argument reduces the problem to the calculation of
$\HF_{n_0}^{[\alpha,\,b)}(H_{\nu})$ for $-2\pi< \alpha <-\pi$ and $\nu$ large
enough. By Lemma \ref{shell}, we have
$$
\HF_{n_0}^{[\alpha,\,b)}(H_{\nu})=H_{n_0+j(1)}(SM,\Z_2)=H_{\dim SM}(SM,\Z_2)=\Z_2
$$
since $n_0+j(1)=\dim SM$. This completes the proof of Claim \ref{claim1}.

\subsubsection{ Proof of Claim \ref{claim2}}

Again by the scaling properties of the functions $H^s_{\nu}$, it suffices 
to prove Claim \ref{claim2}
when $s_0=1$ and $1 < s_1 < \sqrt{1+1/(k'+1)}$. 

We must first find constants $c<a$ for which the homomorphism
\begin{equation}
\labell{eq:c-hat}
\hat{\sigma}_{H^{s_1}_{\nu}H_{\nu}} \colon  \HF^{[a,\,b)}(H_{\nu}) \to
\HF^{[c,\,b)}(H^{s_1}_{\nu})
\end{equation}
is well-defined. For any $c <a$, the function $H^{s_1}_{\nu}$ is in $\Hh^{c,b}$ for all sufficiently large $\nu$. This follows again from the fact that the actions $\Aa^+_{\nu,k}(s_1)$ decrease monotonically to $- \pi k s_1^2$ as 
$\nu\to\infty$.  So, it remains for us to find constants $c<a$ 
such that the Floer chain map $\sigma_{H^{s_1}_{\nu}H_{\nu}}$ satisfies
\begin{equation}
\labell{drop}
\sigma_{H^{s_1}_{\nu}H_{\nu}} ( \CF^a(H_{\nu}) ) \subset 
 \CF^c(H^{s_1}_{\nu})
\end{equation}
(see Section \ref{sec:mon-hom}).
 
In fact, condition \eqref{drop}  only makes 
sense for generic perturbations of $H_{\nu}$ and $H^{s_1}_{\nu}$ 
whose closed orbits, with period equal to one and negative action, 
are nondegenerate. Let $\tilde{H}_{\nu}$ and 
$\tilde{H}^{s_1}_{\nu}$ be such perturbations of $H_{\nu}$ and 
$H^{s_1}_{\nu}$, respectively. Then it suffices to find constants $c<a$ for 
which the condition
\begin{equation}
\labell{tilde-drop}
\sigma_{\tilde{H}^{s_1}_{\nu}\tilde{H}_{\nu}} ( \CF^a(\tilde{H}_{\nu}) ) 
\subset  \CF^c(\tilde{H}^{s_1}_{\nu})
\end{equation}
is satisfied whenever $\tilde{H}_{\nu}$ (respectively, $\tilde{H}^{s_1}_{\nu}$)
 is sufficiently $C^{\infty}$-close to $H_{\nu}$ (respectively, 
$H^{s_1}_{\nu}$).

Consider first the case when $a \in (-k'\pi,\, -(k'+1)\pi)$.
For any $\epsilon >0 $, we can choose $\nu$ sufficiently large
and $\|H_{\nu}-\tilde{H}_{\nu}\|_{C^{\infty}}$ sufficiently small,
so that the actions of the periodic orbits for $\tilde{H}_{\nu}$ 
are all within $\epsilon$ of some $\Aa^{\pm}_{\nu,k}(1)$.   
Similarly, we can choose $\nu$ and $\tilde{H}^{s_1}_{\nu}$
so that the actions of the periodic orbits for $\tilde{H}^{s_1}_{\nu}$ 
are  within $\epsilon$ of some $\Aa^{\pm}_{\nu,k}(s_1)$.

Let $x$ be any periodic orbit of $\tilde{H}_{\nu}$ with period equal to one 
and action less than $a$, i.e., $x \in \Pp^a(\tilde{H}_{\nu})$. By choosing 
$\epsilon
< -k'\pi -a$, we have
\begin{equation}
\labell{ineq1} 
\Aa_{\tilde{H}_{\nu}}(x) \leq -(k'+1)\pi + \epsilon.
\end{equation}

Now, the Floer chain map $\sigma_{\tilde{H}^{s_1}_{\nu}\tilde{H}_{\nu}}$ is 
defined using a homotopy $\tilde{H}^{s}_{\nu}$ from $\tilde{H}_{\nu}$ to 
$\tilde{H}^{s_1}_{\nu}$. We can choose this homotopy so that $ \|\partial_s \tilde{H}^{s}_{\nu}- \partial_s H^{s}_{\nu}\|_{L^{\infty}}$ is arbitrarily small. 

Assume that $y \in \Pp^a(\tilde{H}^{s_1}_{\nu})$ appears in $\sigma_{\tilde{H}^{s_1}_{\nu}\tilde{H}_{\nu}}(x)$ with nonzero coefficient. By the definition of $\sigma_{\tilde{H}^{s_1}_{\nu}\tilde{H}_{\nu}}$, this implies the existence of a perturbed pseudo-holomorphic cylinder $u \colon \R \times S^1 \to W$ whose ends are mapped to $x$ and $y$. The energy of this cylinder satisfies 
\begin{equation*}
\labell{energy}
E(u) \leq  \Aa_{\tilde{H}_{\nu}}(x) - \Aa_{\tilde{H}^{s_1}_{\nu}}(y) + \sup \{\partial_s \tilde{H}^{s}_{\nu}\}.
\end{equation*}
(See, for example,  \cite{bps} Section 4.4.)

Since $ H^{s}_{\nu}$ is monotone and $E(u) \geq 0$, this implies that
\begin{equation}
\labell{ineq2} 
\Aa_{\tilde{H}^{s_1}_{\nu}}(y) \leq \Aa_{\tilde{H}_{\nu}}(x) + \|\partial_s \tilde{H}^{s}_{\nu}- \partial_s H^{s}_{\nu}\|_{L^{\infty}}.
\end{equation}

Together, inequalities \eqref{ineq1} and \eqref{ineq2} yield
\begin{equation}
\labell{ineq3}
\Aa_{\tilde{H}^{s_1}_{\nu}}(y) \leq -(k'+1)\pi + \epsilon + \|\partial_s \tilde{H}^{s}_{\nu}- \partial_s H^{s}_{\nu}\|_{L^{\infty}}.
\end{equation}
Now, the assumption that $ s_1 < \sqrt{1+1/(k'+1)}$ implies that
$$
-(k'+1)\pi < -k'\pi s_1^2.
$$ 
If we choose $\epsilon$ and $ \|\partial_s \tilde{H}^{s}_{\nu}- \partial_s H^{s}_{\nu}\|_{L^{\infty}}$ to be small enough to satisfy
$$
2\epsilon +  \|\partial_s \tilde{H}^{s}_{\nu}- \partial_s H^{s}_{\nu}\|_{L^{\infty}} <  -k'\pi s_1^2 + (k'+1)\pi,
$$
then we get
$$
\Aa_{\tilde{H}^{s_1}_{\nu}}(y) < -k'\pi s_1^2 - \epsilon.
$$
This implies that
$$
\Aa_{\tilde{H}^{s_1}_{\nu}}(y) \leq  -(k'+1)\pi s_1^2 + \epsilon.
$$

Hence,  for $a \in (-k'\pi,\, -(k'+1)\pi)$, condition 
\eqref{tilde-drop} holds for all $c \in (a,\, -(k'+1)\pi s_1^2)$.
 
Next we consider the case when $a= -(k'+1)\pi$. If $\nu$ is sufficiently large
and $\|H_{\nu}-\tilde{H}_{\nu}\|_{C^{\infty}}$ is sufficiently small,
then any periodic orbit $x \in  
CF^a(\tilde{H}_{\nu})$ will satisfy
\begin{equation*}
\labell{ineq4} 
\Aa_{\tilde{H}_{\nu}}(x) \leq -(k'+2)\pi + \epsilon.
\end{equation*}
This follows from the fact that the actions $\Aa^+_{\nu, k}$ decrease
monotonically to $-k\pi$. Using the same arguments as above, 
we see that  any $y \in CF^a(\tilde{H}^{s_1}_{\nu})$ which appears in  $\sigma_{\tilde{H}^{s_1}_{\nu}\tilde{H}_{\nu}}(x)$ with nonzero coefficient satisfies 
\begin{equation*}
\labell{ineq5}
\Aa_{\tilde{H}^{s_1}_{\nu}}(y) \leq -(k'+2)\pi +  \|\partial_s \tilde{H}^{s}_{\nu}- \partial_s H^{s}_{\nu}\|_{L^{\infty}}.
\end{equation*}
The assumption that $ s_1 < \sqrt{1+1/(k'+1)}$ is equivalent to
$$
-(k'+2)\pi < -(k'_1+1)\pi s_1^2,
$$ 
so choosing $\epsilon$ and $ \|\partial_s \tilde{H}^{s}_{\nu}- \partial_s H^{s}_{\nu}\|_{L^{\infty}}$ to be
small enough we get
$$
\Aa_{\tilde{H}^{s_1}_{\nu}}(y) \leq  -(k'+2)\pi s_1^2 + \epsilon.
$$
Thus, when $a= -(k'+1)\pi$, condition \eqref{tilde-drop} holds for 
any $c \in (-(k'+1)\pi,\, -(k'+2)\pi s_1^2)$.

To complete the proof of Claim \ref{claim2} we must prove that 
$\hat{\sigma}_{H^{s_1}_{\nu}H_{\nu}}$ is an isomorphism. This will 
follow almost immediately from Lemma \ref{hat} and the fact that for 
large $\nu$ the actions $ \Aa^+_{\nu,k}(s) $ are arbitrarily
close to $-k \pi s^2$ for all $s \in [1,\, s_1]$. 

For $a \in (-k'\pi,\, -(k'+1)\pi)$, we choose $c \in (a,\, -(k'+1)\pi s_1^2)$ so that $\hat{\sigma}_{H^{s_1}_{\nu}H_{\nu}}$ is well-defined. (Note that here we can choose $c$ so that $c \notin -\pi s_1^2 \Z^+$.) Then, for $\nu$
sufficiently large, it is easy to see that there are continuous families $a_s$ with 
$a_1 =a$ and $a_{s_1}=c$, such that 
$$
\Aa^+_{\nu,k'+1}(s) <a_s  <\Aa^+_{\nu,k'}(s)
$$ 
for all $s \in [1,\, s_1]$. For any such family $a_s$, $\tilde{H}^s_{\nu} \in \Hh^{a_s,b}$ for all $s \in [1,\, s_1]$. Lemma \ref{hat} then implies that 
$\hat{\sigma}_{H^{s_1}_{\nu}H_{\nu}}$ is an isomorphism. A similar argument 
works for the case $a= -(k'+1)\pi$.

This completes the proof of Claim \ref{claim2} and hence Proposition \ref{cont}.

\section{Homological symplectic capacity}
\labell{sec:capacity}
Following \cite{fhw} closely, we will now define a certain relative 
homological capacity. The results of of this paper imply the finiteness 
of this capacity in various cases.

As above, let $W$ be a geometrically bounded symplectically aspherical
manifold. For an open subset $U\subset W$, we set
$$
\SH^{[a,\,0)}(U)= \varinjlim_{b\nearrow 0}\SH^{[a,\,b)}(U).
$$
Furthermore, for a compact subset $Z\subset W$, we define its symplectic 
homology as
$$
\SH^{[a,\,0)}(Z)= \varprojlim_{U\searrow Z}\SH^{[a,\,0)}(U),
$$
where the limit is taken over all open sets $U\supset Z$. By the definition
of the inverse limit, we have the projection map
$$
\phi^a_U\colon \SH^{[a,\,0)}(Z)\to \SH^{[a,\,0)}(U)
$$
for every open set $U\supset Z$. 

\begin{Example}
Let $M$ be a closed symplectic submanifold of $W$ and let $B_R$ be the 
tubular neighborhoods of $M$ constructed in Section
\ref{sec:pf}. As readily follows from
Proposition \ref{cont}, $\SH^{[a,\,0)}_{n_0}(M)=\Z_2$, where
$n_0=\frac{3}{2}\dim W-\codim M -1$. Moreover,
$\phi^a_{B_R}\colon \SH^{[a,\,0)}_{n_0}(M)\to \SH^{[a,\,0)}_{n_0}(B_R)$ 
is an isomorphism, provided that $a<-\pi R^2$. Furthermore,
$\SH^{[a,\,0)}(B_R)=0$ if $-\pi R^2 <a<0$.  
\end{Example}

Define the \emph{relative homological capacity} of $(U,Z)$ as
$$
\wsh (U,Z)=\inf\{-a\mid \phi^c_U\neq 0\text{ for all } c<a\}.
$$
Clearly, $\wsh$ is an invariant of symplectomorphisms of the ambient manifold
$W$. The fact that $\wsh$ is a symplectic capacity (defined only on 
submanifolds of $W$) can be easily verified using Proposition \ref{cont}. In
other words, we have

\begin{Theorem}
\labell{thm:wsh}
~
\begin{enumerate}

\item {\rm [Invariance].} The relative capacity $\wsh$ is an invariant of 
symplectomorphisms of $W$.

\item {\rm [Monotonicity].} Let $Z'\subset Z\subset U\subset U'$. Then
$\wsh(U,Z)\leq \wsh(U',Z)$ and $\wsh(U,Z)\leq \wsh(U,Z')$.
\item {\rm[Homogeneity].} For any constant $a>0$,
$$
\wsh(U,Z,a\omega)=a \wsh(U,Z,\omega).
$$
\item {\rm [Normalization].} Assume that $M$ is a closed symplectic 
submanifold of a geometrically bounded symplectically aspherical manifold 
$W$ and let $B_R$ be a symplectic tubular neighborhood of $M$ 
in $W$ of radius $R>0$ (see Section \ref{sec:pf}).
Then $\wsh(B_R,M)=\pi R^2$.
\end{enumerate}
\end{Theorem}

Here the first three assertions are obvious. To prove the
last assertion note that $\wsh(B_R,M)\leq \pi R^2$ by Proposition \ref{cont}. 
On the other hand, $\SH^{[a,\, 0)}(B_R)=0$ when $-\pi R^2 < a<0$ as readily
follows from the analysis carried out in Section \ref{sec:exh} (see
Lemma \ref{empty}). Thus,
$\wsh(B_R,M)\geq \pi R^2$, which completes the proof of the last assertion.

The capacity $\wsh(U,\point)$ is essentially the homological capacity
introduced in \cite{fhw}. Once, $\wsh(U,Z)<\infty$,
the symplectic homology of $V$ is non-zero for any open subset $V$ of $U$
such that $Z\subset V\subset U$. Then a suitable version of
Proposition \ref{near}
implies the nearby existence theorem in $(U,Z)$, i.e., the existence of
contractible periodic orbits on a dense of levels of any function
attaining an isolated minimum on $Z$. Yet, the 
finiteness of this capacity does not seem to imply the existence of 
periodic orbits for almost all energy levels as does the Hofer--Zehnder 
capacity.
A different version of a relative capacity, a relative analogue of the
Hofer--Zehder capacity, which is sufficient for proving a version of 
the almost existence theorem is introduced and analyzed in detail in \cite{gg}.

Finally, we note that the relative homological capacity $\wsh$ is different from, although apparently related to, the one introduced in \cite{bps}.

\end{document}